\documentclass[11pt]{article}

\usepackage{amsmath,amssymb,amsthm}
\usepackage{booktabs}
\usepackage{graphicx}
\usepackage{microtype}
\usepackage{geometry}
\usepackage[authoryear,round]{natbib}
\usepackage[colorlinks=true,linkcolor=blue,citecolor=blue,urlcolor=blue]{hyperref}

\usepackage{setspace}

\usepackage{doi}

\newcommand{\E}{\mathbb{E}}
\newcommand{\Var}{\mathrm{Var}}
\newcommand{\Cov}{\mathrm{Cov}}
\newcommand{\GMD}{\mathrm{GMD}}

\newtheorem{proposition}{Proposition}
\newtheorem{corollary}{Corollary}

\title{Order-Induced Variance in the Moving-Range Sigma Estimator: A Total-Variance Decomposition}
\author{Andrew T. Karl\\Karl Statistical Services LLC\\\texttt{akarl@asu.edu}}
\date{}

\begin{document}
\maketitle

\vspace{1em}
\begin{abstract}
Individuals and Moving Range (I--MR) charts commonly estimate the process standard deviation $\sigma$ by dividing the span-2 average moving range by the Normal-reference unbiasing constant $d_2$.  Unlike the bias-corrected sample standard deviation $S/c_4$, this estimator depends on ordering through adjacency, so permuting a fixed sample changes the estimate.  We formalize this dependence by introducing an independent uniformly random permutation and applying the law of total variance.  This gives an exact decomposition of the sampling variance into a values component (the variance of the permutation mean) and an adjacency component (the expected conditional variance over permutations).  The permutation mean is order-invariant and equals $\GMD/d_2$, where $\GMD$ is the sample Gini mean difference.  Under i.i.d. Normal sampling, both components admit closed forms.  Under the Normal reference, this adjacency component accounts for nearly all of the familiar asymptotic variance inflation of MR(2)/$d_2$ relative to $S/c_4$.  The decomposition itself is not Normal-specific: an appendix gives exact finite-sample expressions for any i.i.d. distribution with finite variance.  
\end{abstract}

\begin{center}
\begin{minipage}{0.9\textwidth}
\small
\textbf{Publication notice.}
This manuscript was accepted for publication and has since been published
in its final form as:

\medskip

\noindent
Karl, A. T. (2026).
``Order-induced variance in the moving-range sigma estimator:
a total-variance decomposition.''
\emph{Statistical Papers}, \textbf{67}, Article 97.
\href{https://doi.org/10.1007/s00362-026-01877-0}
{https://doi.org/10.1007/s00362-026-01877-0}.

\medskip

\noindent
The published version of record may differ from this author manuscript.
\end{minipage}
\end{center}

\noindent\textbf{Keywords:} statistical process control; I--MR charts; permutation distribution; law of total variance; Gini mean difference; asymptotic relative efficiency

\section{Motivation}
\label{sec:motivation}

The standard Individuals and Moving Range (I--MR) estimate of the process standard deviation $\sigma$ is the span-2 average moving range divided by the usual constant $d_2$ \citep{Montgomery2019ISQC}. Under i.i.d. Normal sampling, $d_2=\E|X_1-X_2|/\sigma=2/\sqrt{\pi}$, so the span-2 moving-range estimator MR(2)/$d_2$ is unbiased for $\sigma$ but less efficient than $S/c_4$ \citep{Hoel1946,CryerRyan1990,BraunPark2008}. Because MR(2) uses only successive differences, it depends on ordering through adjacency. In practice, MR-based summaries are used as short-term dispersion measures and are sometimes preferred in Phase~I when drift or oscillation inflates $S$ \citep[e.g.,][]{Nelson1982,CryerRyan1990}; see \citet{RigdonCruthisChamp1994} for design strategies for I--MR charts.

The SPC literature on dispersion monitoring includes exponentially weighted moving average (EWMA) and other memory-type charts for process dispersion \citep{HuwangHuangWang2010,AliHaq2017}, nonparametric sign-based dispersion charts \citep{GodaseEtAl2022}, lognormal-process variability monitoring with estimated parameters \citep{AkhtarEtAl2022QRE,AkhtarEtAl2024QRE}, ranked-set-sampling dispersion charts \citep{AbbasEtAl2023}, and multivariate variability or covariance monitoring \citep{AjadiEtAl2021QRE,AjadiEtAl2024QRE}. These contributions primarily address chart design, run-length performance, and robustness of monitoring schemes; \citet{JensenEtAl2006} review the effects of parameter estimation on control-chart properties. The aim here is narrower: the sampling variance of the classical span-2 moving-range estimator and the portion of that variance induced by using only adjacent pairs.

With the observed values held fixed, different permutations of the same sample generally produce different average moving ranges. \citet{PootsWoodcock2012} give the attainable minimum and maximum of the average moving range over all permutations. In time-ordered applications this sensitivity is intentional, but it still contributes to sampling variability even when the data are i.i.d. In practical terms, the adjacency component is the precision cost of replacing all pairwise differences by the $n-1$ local differences along one path through the data. If the observed order is arbitrary or effectively random, this component is pure order-induced noise. If the order is meaningful, as in a time sequence, the same component is also the part of the estimator that can respond to local smoothness, roughness, or serial dependence. Thus the decomposition is not a prescription to ignore time order; it separates marginal scale information in the values from the additional contribution of the particular ordering.

This ``values versus order'' distinction is explicit in Shewhart's original framing of control-chart evidence \citep[p.~90]{Shewhart1939}. He emphasizes that the informational content of a sequence is not exhausted by its marginal distribution and that order can be a primary diagnostic clue \citep[pp.~25--27]{Shewhart1939}. To make this concrete, he considers exactly the thought experiment relevant here: take a fixed set of observed numbers, imagine the order being lost by mixing the items, and compare the observed production order against the ensemble of random orders generated by repeated mixing and drawing \citep[pp.~31--32]{Shewhart1939}. In modern terms, we condition on $X=x$ and randomize the permutation. In his Fig.~7 example, he notes that most of the possible reorderings would fail to reveal the assignable-cause signal that is visible in the production order: ``it was the order that furnished the clue'' \citep[pp.~31--33, 90]{Shewhart1939}.

We therefore introduce an independent uniformly random permutation $\Pi$ and view the moving-range estimator as $T(X,\Pi)$. The law of total variance decomposes $\Var\{T(X,\Pi)\}$ into a values component, $\Var\{\E(T\mid X)\}$, and an adjacency component, $\E\{\Var(T\mid X)\}$. The permutation mean is order-invariant and equals the sample Gini mean difference divided by $d_2$. Under i.i.d. Normal sampling, both components are available in closed form, allowing a direct interpretation of how much of $\Var(T)$ is attributable to adjacency; exact distribution-free expressions for the components are collected in Appendix~\ref{app:exact}.

\section{Setup: the moving-range functional with random ordering}
\label{sec:setup}

Let $X_1,\dots,X_n$ $(n\ge 2)$ be i.i.d. with finite common standard deviation $\sigma$; write $X=(X_1,\dots,X_n)$, and let $\Pi$ be an independent uniformly random permutation of $\{1,\dots,n\}$. Define the unscaled average moving range (the MR(2) statistic of Section~\ref{sec:motivation}) along order $\Pi$ by
\[
\overline{\mathrm{MR}}(X,\Pi)
:=\frac{1}{n-1}\sum_{t=2}^{n}\bigl|X_{\Pi(t)}-X_{\Pi(t-1)}\bigr|.
\]
The decomposition in Section~\ref{sec:lotv} is distribution-free under i.i.d. sampling with finite variance. Normality enters only through the Normal-reference unbiasing constants ($d_2$, $c_4$) and the closed-form evaluations in Section~\ref{sec:normal}.

Define the scaled moving-range functional
\begin{equation}
T(X,\Pi)=\frac{\overline{\mathrm{MR}}(X,\Pi)}{d_2}.
\label{eq:Tdef}
\end{equation}
For a fixed observed ordering $\Pi_{\mathrm{obs}}$ (e.g., time order), the usual I--MR estimate is $T_{\mathrm{obs}}:=T(X,\Pi_{\mathrm{obs}})$. Throughout, $\Pi$ is purely a mathematical device for separating the contribution of the observed values from that of their arrangement; we neither assume that the observed ordering arose by random permutation nor propose discarding time order in applications (Section~\ref{sec:discussion} returns to this point). Write $\E_{\Pi}(\cdot\mid X)$ and $\Var_{\Pi}(\cdot\mid X)$ for expectation and variance over $\Pi$ holding $X$ fixed, and define the permutation mean
\[
\bar T(X):=\E_{\Pi}\{T(X,\Pi)\mid X\}.
\]

\section{A law of total variance decomposition}
\label{sec:lotv}
The statistic $T(X,\Pi)$ varies due to both the realized values $X$ and, conditional on $X$, the adjacency pattern induced by the ordering $\Pi$. The law of total variance separates these contributions exactly.

\subsection{Values and adjacency components}
\label{subsec:lotv_split}
\begin{proposition}[Total-variance decomposition]
For $(X,\Pi)$ and $T(X,\Pi)$ as in Section~\ref{sec:setup},
\begin{equation}
\Var\bigl\{T(X,\Pi)\bigr\}
=
\E\Bigl[\Var\bigl\{T(X,\Pi)\mid X\bigr\}\Bigr]
+
\Var\Bigl[\E\bigl\{T(X,\Pi)\mid X\bigr\}\Bigr].
\label{eq:lotv}
\end{equation}
\end{proposition}
\noindent This is the law of total variance applied to $T(X,\Pi)$ (see, e.g., \citealp[Thm.~4.4.7]{CasellaBerger2002}); the variance is finite because each summand in $\overline{\mathrm{MR}}$ has variance at most $\E(X_1-X_2)^2=2\sigma^2$.

The first term is an adjacency component: it is the average, over the distribution of $X$, of the fixed-sample permutation variance $\Var_{\Pi}\{T(x,\Pi)\}$. This term is not ``process variance'' in the usual SPC sense; it is the expected conditional variance induced by randomizing the order while holding the observed values fixed. The second term is a values component, $\Var\{\bar T(X)\}$: it reflects variation of the order-invariant permutation mean $\bar T(X)$ across realizations of $X$ and is therefore determined by the values alone, not their arrangement.

Because $(X_1,\dots,X_n)$ is exchangeable under i.i.d. sampling, $T(X,\Pi_{\mathrm{obs}})\overset{d}=T(X,\Pi)$ for any fixed ordering $\Pi_{\mathrm{obs}}$; thus $\Var\{T(X,\Pi)\}$ in \eqref{eq:lotv} equals the sampling variance of the conventional MR(2) estimator.

A convenient summary is the adjacency fraction
\begin{equation}
\mathrm{AdjFrac}(n)
:=
\frac{\E\bigl[\Var\{T(X,\Pi)\mid X\}\bigr]}{\Var\{T(X,\Pi)\}}.
\label{eq:adjfrac}
\end{equation}
The fraction depends on the sampling distribution as well as on $n$; Section~\ref{sec:normal} evaluates it under the Normal reference, and Corollary~\ref{cor:adjfrac} shows that the dependence on the distribution is only through the correlation $\rho_1$ defined in Appendix~\ref{app:exact}. Since $T=\overline{\mathrm{MR}}/d_2$ is a constant rescaling, $\mathrm{AdjFrac}(n)$ is unchanged if defined using $\overline{\mathrm{MR}}$ instead of $T$.

\subsection{Permutation mean as a Gini baseline}
\label{subsec:gmd}
\begin{proposition}[Permutation mean and Gini mean difference]\label{prop:gmd}
For any fixed realization $x=(x_1,\dots,x_n)$,
\begin{equation}
\bar T(x)
=\frac{\GMD(x)}{d_2},
\label{eq:gmd}
\end{equation}
where $\GMD(x)=\frac{2}{n(n-1)}\sum_{i<j}|x_i-x_j|$ is the sample Gini mean difference \citep{David1968,GerstenbergerVogel2015}.
\end{proposition}
\noindent The proof is given in Appendix~\ref{app:prop2}.

Equation~\eqref{eq:gmd} shows that averaging over adjacencies replaces the average over the $n-1$ adjacent pairs in $\overline{\mathrm{MR}}$ with the average over all $\binom{n}{2}$ pairs. Thus the Gini mean difference is the natural random-order baseline for the moving range. Equivalently, $\bar T$ is a symmetric degree-two $U$-statistic up to the fixed factor $1/d_2$ \citep{David1968,GerstenbergerVogel2015}. Gini-based dispersion measures have also appeared directly in control-chart design; see, for example, \citet{RiazSaghir2007} and \citet{AslamEtAl2022}. More broadly, pairwise-difference scale functionals, including Gini's mean difference and generalized $Q_n$ statistics, have been used for robust tests of scale changes in time series \citep{GerstenbergerVogelWendler2020}.

\section{Normal-reference expressions and an efficiency interpretation}
\label{sec:normal}
Under i.i.d. $N(\mu,\sigma^2)$ sampling, $d_2=2/\sqrt{\pi}$ (so $1/d_2^2=\pi/4$), and both terms in \eqref{eq:lotv} admit closed forms.

\subsection{Closed-form component variances}
\label{subsec:closedform}
\paragraph{Values component.}
Using the Normal closed forms for $\Var\{\GMD(X)\}$ \citep{ZengaPolisicchioGreselin2004,GerstenbergerVogel2015},
\begin{equation}
\Var\{\bar T(X)\}
=
\frac{\pi}{4}\,\Var\{\GMD(X)\}
=
\frac{\sigma^2}{n(n-1)}
\left[
\frac{\pi(n+1)}{3}
+2\sqrt{3}(n-2)-2(2n-3)
\right].
\label{eq:VarTbar_normal_sigma}
\end{equation}

\paragraph{Total variance.}
As noted after \eqref{eq:lotv}, $\Var\{T(X,\Pi)\}$ equals the variance of $T(X,\Pi_{\mathrm{obs}})$ along any fixed ordering $\Pi_{\mathrm{obs}}$, in particular the observed-order estimator $T_{\mathrm{obs}}$.
Let
\[
 w \,=\, \overline{\mathrm{MR}}(X,\Pi_{\mathrm{obs}})=\frac{1}{n-1}\sum_{t=2}^{n}|X_t-X_{t-1}|,
\]
taking $\Pi_{\mathrm{obs}}$ to be the identity without loss of generality; this matches the definition of \citet[eq.~(3)]{Hoel1946}. Under i.i.d. Normal sampling, Hoel's variance formula \citep[eq.~(10)]{Hoel1946}, multiplied by $\sigma^2$ for general scale and by $1/d_2^2=\pi/4$, gives
\begin{equation}
\Var\{T(X,\Pi)\}
=
\frac{\pi}{4}\,\Var(w)
=
\frac{\pi\,\sigma^2}{2\,(n-1)^2}\left[
\left(\frac{4}{3}+\frac{2\sqrt{3}-6}{\pi}\right)n
+\left(\frac{10-4\sqrt{3}}{\pi}-\frac{5}{3}\right)
\right].
\label{eq:VarT_normal_sigma}
\end{equation}
The adjacency component is the difference \eqref{eq:VarT_normal_sigma}\,$-$\,\eqref{eq:VarTbar_normal_sigma}; Table~\ref{tab:components} reports all three quantities.

\subsection{Adjacency fraction and an efficiency interpretation}
\label{subsec:eff}
From \eqref{eq:adjfrac} and \eqref{eq:lotv},
\begin{equation}
\mathrm{AdjFrac}(n)=1-\frac{\Var\{\bar T(X)\}}{\Var\{T(X,\Pi)\}}.
\label{eq:adjfrac_equiv}
\end{equation}
Let $a=\pi/3+2\sqrt3-4$ and $c=4/3+(2\sqrt3-6)/\pi$, so that, as $n\to\infty$, $\Var\{\bar T(X)\}\sim \sigma^2 a/n$ and $\Var\{T(X,\Pi)\}\sim \sigma^2(\pi c/2)/n$.  Hence
\begin{equation}
\mathrm{AdjFrac}(\infty):=\lim_{n\to\infty}\mathrm{AdjFrac}(n)
=
1-\frac{2a}{\pi c}
=
\frac{\pi+3-3\sqrt{3}}{2\pi+3\sqrt{3}-9}
\approx 0.3813.
\label{eq:AdjFrac_limit_normal}
\end{equation}
Thus, even with i.i.d. data, roughly $38\%$ of the moving-range estimator's sampling variance is attributable to random adjacency. Appendix~\ref{app:exact} shows that this limit is the Normal instance of a distribution-free formula: $\mathrm{AdjFrac}(n)$ depends on the sampling distribution only through the correlation $\rho_1$ between two absolute pairwise differences sharing a common observation, with limit $(1-2\rho_1)/(1+2\rho_1)$; the Normal value is $\rho_1\approx 0.2239$.

Hoel \citeyearpar{Hoel1946} reported $\mathrm{ARE}(T,S)\approx 0.605$ under Normality. The present decomposition explains where this loss arises. The standardized Gini mean difference is nearly as efficient as $S$, with $\mathrm{ARE}(\bar T,S)\approx 0.978$ under Normality \citep[Table~6]{GerstenbergerVogel2015}, while
\[
\frac{\Var(\bar T)}{\Var(T)}\to 1-\mathrm{AdjFrac}(\infty)\approx 0.6187.
\]
Because ARE is a ratio of asymptotic variances, the Normal-reference factorization is
\[
\mathrm{ARE}(T,S)
=
\mathrm{ARE}(\bar T,S)\,\bigl(1-\mathrm{AdjFrac}(\infty)\bigr)
\approx 0.978\times 0.6187\approx 0.605.
\]
At finite $n$, the analogous statement is the exact variance-ratio identity obtained from \eqref{eq:adjfrac_equiv} and Table~\ref{tab:components}; the limiting constants $0.978$ and $0.6187$ should not be substituted at small $n$.

Thus, MR(2)/$d_2$ is inefficient not primarily because it uses absolute differences rather than squared deviations---the Gini baseline is nearly as efficient as $S$---but because it uses only $n-1$ adjacent pairs rather than all $\binom n2$ pairs. Equivalently, $1/\mathrm{ARE}(T,S)\approx 1.653$ while $1/\mathrm{ARE}(\bar T,S)\approx 1.023$, so adjacency accounts for about $(1.653-1.023)/(1.653-1)\approx 96.5\%$ of the asymptotic variance inflation of $T$ beyond $S$ (or $S/c_4$).

\subsection{Reference magnitudes}
Table~\ref{tab:components} reports Normal-reference values for $\Var\{T(X,\Pi)\}$ and its decomposition at representative $n$.  
For comparison, the table also includes $\Var\{S/c_4\}=(1-c_4^{2})/c_4^{2}$ \citep{CryerRyan1990}.

\begin{table}[t]
\centering
\caption{Normal-reference variance decomposition for $T=\overline{\mathrm{MR}}/d_2$ under $N(0,1)$ sampling. Columns give the total variance \eqref{eq:VarT_normal_sigma}, the adjacency component \eqref{eq:VarT_normal_sigma}\,$-$\,\eqref{eq:VarTbar_normal_sigma}, the values component \eqref{eq:VarTbar_normal_sigma}, the adjacency fraction \eqref{eq:adjfrac}, and, for comparison, $\Var\{S/c_4\}$.}
\label{tab:components}
\begin{tabular}{cccccc}
\toprule
$n$ & $\Var\{T(X,\Pi)\}$ & $\E\{\Var(T\mid X)\}$ & $\Var\{\bar T(X)\}$ & $\mathrm{AdjFrac}(n)$ & $\Var\{S/c_4\}$ \\
\midrule
4   & 0.2471 & 0.0667 & 0.1803 & 0.2701 & 0.1781 \\
8   & 0.1128 & 0.0377 & 0.0752 & 0.3339 & 0.0738 \\
12  & 0.0730 & 0.0256 & 0.0474 & 0.3511 & 0.0464 \\
16  & 0.0540 & 0.0194 & 0.0346 & 0.3591 & 0.0339 \\
20  & 0.0428 & 0.0156 & 0.0272 & 0.3638 & 0.0267 \\
25  & 0.0340 & 0.0125 & 0.0215 & 0.3674 & 0.0210 \\
50  & 0.0168 & 0.0063 & 0.0105 & 0.3745 & 0.0103 \\
100 & 0.0083 & 0.0031 & 0.0052 & 0.3779 & 0.0051 \\
\bottomrule
\end{tabular}
\end{table}

The limiting value in \eqref{eq:AdjFrac_limit_normal} should be interpreted as a large-sample Normal-reference benchmark rather than as a replacement for the finite-$n$ values in Table~\ref{tab:components}. For typical Phase~I sample sizes the convergence is already fairly close: for example, $\mathrm{AdjFrac}(20)=0.364$, $\mathrm{AdjFrac}(35)=0.371$, and $\mathrm{AdjFrac}(50)=0.374$, compared with the limiting value $0.381$. The consequence---that a substantial share of the sampling variance of MR(2)/$d_2$ is attributable to adjacency---is already visible at common Phase~I sample sizes.

\section{Discussion}
\label{sec:discussion}
The moving-range estimator MR(2)/$d_2$ depends on two ingredients: the realized values and the induced adjacency pattern. Introducing an independent random permutation and applying the law of total variance yields the exact split \eqref{eq:lotv} into values and adjacency components, with permutation mean $\bar T=\GMD/d_2$. Under Normality, the adjacency component is substantial and accounts for nearly all of the familiar asymptotic variance inflation relative to $S/c_4$ (Section~\ref{subsec:eff}). This formalizes Shewhart's long-standing distinction between the numbers observed and their order, and makes precise how much sampling variability arises from adjacency alone \citep[p.~90]{Shewhart1939}.

The decomposition is not a Normal artifact: Appendix~\ref{app:exact} gives exact finite-sample expressions for both components in \eqref{eq:lotv} in terms of $\zeta_2=\Var|X_1-X_2|$ and $\zeta_1=\Cov(|X_1-X_2|,|X_1-X_3|)$, and by Corollary~\ref{cor:adjfrac} the adjacency fraction depends on the sampling distribution only through $\rho_1=\zeta_1/\zeta_2\in[0,1/2]$, with the limit in \eqref{eq:adjfrac_rho}. Representative limiting values are $2/3$ for the Uniform, $0.3813$ for the Normal, approximately $0.23$ for the heavy-tailed $t_5$, and $1/5$ for the skewed Exponential. Thus a substantial adjacency share is not confined to the Normal case, although its magnitude is distribution-specific. In these examples the heavier-tailed and skewed cases have larger $\rho_1$ and hence smaller adjacency shares (but no monotone tail-weight claim is intended). The short-tailed case is especially instructive: at the Uniform the standardized Gini baseline is exactly as efficient as $S$ \citep[Table~6]{GerstenbergerVogel2015}, so the factorization of Section~\ref{subsec:eff}, with the moving range calibrated by the Uniform value of $\E|X_1-X_2|/\sigma$, gives a limiting $\mathrm{ARE}$ of $1\times\tfrac13=\tfrac13$ versus $S$---a larger efficiency loss than the Normal case ($0.605$), and one attributable entirely to adjacency. For stable processes with bounded, short-tailed variation, the Normal benchmark $0.3813$ can therefore understate the adjacency share, as the Uniform case illustrates.

The fixed divisor $d_2=2/\sqrt{\pi}$ is retained throughout because the target application is the ordinary Normal-reference I--MR estimator. Any fixed divisor cancels from $\mathrm{AdjFrac}(n)$; the separate question of calibrating an unbiased scale estimator under a non-Normal law involves the distribution-specific constant $\E|X_1-X_2|/\sigma$ and is outside the scope of this note. Consequently, the value $0.3813$ and the efficiency comparison with $S/c_4$ should not be exported mechanically to skewed or heavy-tailed settings. In such settings the Gini baseline may itself have different efficiency behavior relative to $S$ \citep{GerstenbergerVogel2015}.

In many applications the order dependence of MR(2)/$d_2$ is intentional, because the moving range targets local variation in a time-ordered sequence. The present results are therefore not an argument against MR charts; rather, they quantify an intrinsic precision cost of localization that persists even when successive observations are effectively i.i.d. The unscaled MR(2) functional is the mean absolute successive difference, paralleling the mean square successive difference of \citet{vonNeumann1941}; Normal-sampling moments of the mean absolute successive difference, and of its ratio to the root mean square, were obtained by \citet{Kamat1953}. Where both $\overline{\mathrm{MR}}/d_2$ and $S/c_4$ are plausible estimators for individuals data, the decomposition quantifies the relevant tradeoff: whether the diagnostic value of adjacency is worth the extra sampling variation from estimating scale along a single adjacent path through the observations. A practical implication is that $\GMD/d_2$ provides a random-order baseline for the observed moving-range estimate. Large discrepancies among $T_{\mathrm{obs}}$, $\bar T$, and $S/c_4$ can then be read as evidence that ordering, local dependence, or nonstationarity is materially affecting the moving-range estimate.

\begin{figure}[t] \centering \includegraphics[width=0.84\linewidth]{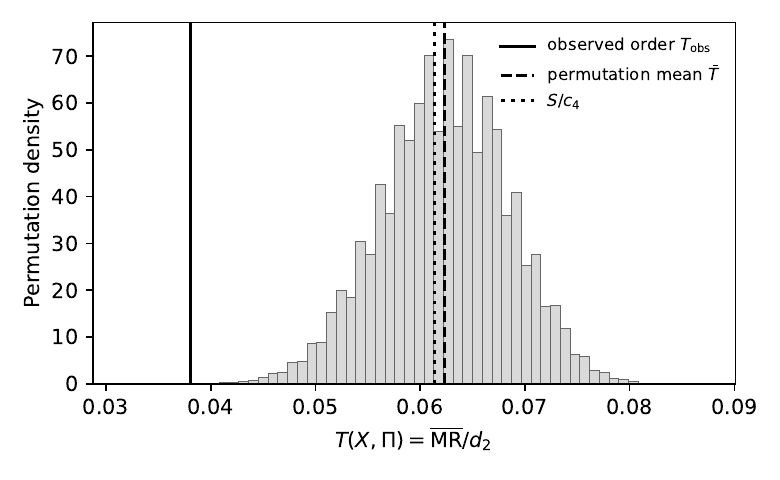} \caption{Conditional permutation distribution of $T(X,\Pi)=\overline{\mathrm{MR}}/d_2$ for the $n=35$ chemical-process data of \citet[Table~2]{CryerRyan1990}, based on $B=10^{7}$ random permutations of the observed values. The solid vertical line marks the observed time-order estimate $T_{\mathrm{obs}}$, the dashed line the permutation mean $\bar T=\GMD/d_2$, and the dotted line $S/c_4$; for these data $\bar T$ and $S/c_4$ nearly coincide.} \label{fig:cr-permutation} \end{figure}
Because $T(X,\Pi)$ is defined for any permutation, the conditional permutation distribution of $T$ given the observed values gives a descriptive benchmark for comparing the observed order with random reorderings \citep{PootsWoodcock2012}. This benchmark turns the mixing-and-drawing thought experiment of Section~\ref{sec:motivation} into a computable diagnostic \citep[pp.~31--33]{Shewhart1939}. It should not be interpreted as a definitive i.i.d. test from a single realized run; rather, a one-sided summary such as
\[
p=\Pr_{\Pi}\{T(X,\Pi)\le T_{\mathrm{obs}}\mid X\}
\]
can be reported as a diagnostic of whether the observed ordering is unusually smooth or rough relative to random reorderings. The squared analogue has classical precedent: \citet{Young1941} derived the permutation moments of a mean-square-successive-difference ratio for a fixed sample and applied the resulting randomness test to quality-control data.

As an illustration, consider the $n=35$ chemical-process observations in \citet[Table~2]{CryerRyan1990}, where consecutive measurements have sample lag-1 autocorrelation $0.528$. The observed-order estimate is $T_{\mathrm{obs}}=0.038$ (as reported by Cryer and Ryan), while the permutation mean is $\bar T=\GMD/d_2=0.062$, close to $S/c_4=0.061$. In $B=10^{7}$ random reorderings of the same values, the observed estimate lies far in the lower tail of the conditional permutation distribution, with estimated one-sided permutation probability approximately $4\times 10^{-5}$ and Monte Carlo standard error about $2\times 10^{-6}$; see Fig.~\ref{fig:cr-permutation}. Positive serial dependence makes neighbors unusually similar, suppressing $T_{\mathrm{obs}}$ far below the Gini baseline. This contrast complements the recommendation of \citet{CryerRyan1990} to compute both $\overline{\mathrm{MR}}/d_2$ and $S/c_4$ and to investigate substantial discrepancies; see also \citet{CruthisRigdon1992}. In current practice this comparison is summarized by the stability ratio of \citet{RamirezRunger2006} and the stability index of \citet{JensenSzarkaWhite2019}.

\section*{Statements and Declarations}

\noindent\textbf{Funding.}
No funding was received to assist with the preparation of this manuscript.

\medskip
\noindent\textbf{Competing interests.}
The author reports there are no competing interests to declare.

\medskip
\noindent\textbf{Data and code availability.}
The numerical illustration uses the 35 observations reported by \citet[Table~2]{CryerRyan1990}. Online Resource~1 provides reproducible code for the closed-form evaluations, Table~\ref{tab:components}, and the illustration; Online Resource~2 provides a nested-permutation check of the Appendix~\ref{app:exact} component formulas. Both online resources are archived in Mendeley Data \citep{KarlOR2026}.

\medskip
\noindent\textbf{Use of generative AI and AI-assisted technologies.}
During the preparation and revision of this work, the author used OpenAI's ChatGPT 5.5 Pro and Anthropic's Claude Fable 5 to assist with proofreading and rewording text, checking mathematical derivations and numerical results, refining LaTeX formatting, and structuring and commenting the supplementary code, based on the author's instructions. The author reviewed and edited all AI-assisted output and takes full responsibility for the content of the article.

\appendix
\renewcommand{\theequation}{\thesection.\arabic{equation}}
\renewcommand{\theHequation}{\thesection.\arabic{equation}}
\setcounter{equation}{0}
\section{Proof of the permutation-mean identity}
\label{app:prop2}
Fix $x=(x_1,\dots,x_n)$. By definition \eqref{eq:Tdef},
\[
\bar T(x)=\E_{\Pi}\{T(x,\Pi)\}
=\frac{1}{d_2}\cdot\frac{1}{n-1}\sum_{t=2}^{n}
\E_{\Pi}\bigl[\,|x_{\Pi(t)}-x_{\Pi(t-1)}|\,\bigr].
\]
Fix $t\in\{2,\dots,n\}$. Since $\Pi$ is a permutation, $\Pi(t-1)\neq \Pi(t)$.
For any distinct indices $i\neq j$, the number of permutations $\pi$ of
$\{1,\dots,n\}$ satisfying $\pi(t-1)=i$ and $\pi(t)=j$ is $(n-2)!$.
Thus
\[
\Pr\bigl((\Pi(t-1),\Pi(t))=(i,j)\bigr)=\frac{(n-2)!}{n!}=\frac{1}{n(n-1)}
\qquad(i\neq j).
\]
Therefore
\[
\E_{\Pi}\bigl[\,|x_{\Pi(t)}-x_{\Pi(t-1)}|\,\bigr]
=\sum_{i\neq j}|x_i-x_j|\Pr\bigl((\Pi(t-1),\Pi(t))=(i,j)\bigr)
=\frac{1}{n(n-1)}\sum_{i\neq j}|x_i-x_j|.
\]
The right-hand side does not depend on $t$, so the factor $1/(n-1)$ cancels
against the $n-1$ identical summands, yielding
\begin{equation*}
\bar T(x)=\frac{1}{d_2}\cdot\frac{1}{n(n-1)}\sum_{i\neq j}|x_i-x_j|
=\frac{1}{d_2}\cdot\frac{2}{n(n-1)}\sum_{i<j}|x_i-x_j|
=\frac{\GMD(x)}{d_2}.\tag*{$\square$}
\end{equation*}

\section{Exact component variances under general i.i.d. sampling}
\label{app:exact}
\setcounter{equation}{0}

This appendix gives the finite-sample component variances corresponding to \eqref{eq:lotv} without assuming Normality. The formulas are stated for the unscaled moving-range average $\overline{\mathrm{MR}}$; this convention keeps the distributional quantities separate from the Normal-reference constant used in ordinary I--MR practice.

\subsection*{Notation}
Write $D_{ij}:=|X_i-X_j|$ and define
\[
\zeta_2:=\Var(D_{12}),
\qquad
\zeta_1:=\Cov(D_{12},D_{13}),
\qquad
\rho_1:=\zeta_1/\zeta_2=\mathrm{Corr}(D_{12},D_{13})\qquad(\zeta_2>0).
\]
Thus $\zeta_2$ measures the marginal variability of one absolute pairwise
difference, while $\zeta_1$ measures the dependence induced when two such
differences share one observation. By exchangeability, $\zeta_1$ also equals
$\Cov(D_{12},D_{23})$, the lag-one autocovariance of the successive absolute
differences entering $\overline{\mathrm{MR}}$. Moreover,
\[
\zeta_1=\Var\{\eta(X_1)\}\ge 0,
\qquad
\eta(x):=\E|x-X_2|.
\]
Because $D_{12}$ and $D_{13}$ are conditionally independent given $X_1$
(so their conditional covariance vanishes) and have common conditional mean
$\eta(X_1)$, the covariance $\zeta_1$ equals $\Var\{\eta(X_1)\}$, the first
projection variance for the degree-two $U$-statistic with kernel
$h(x,y)=|x-y|$ \citep{Hoeffding1948,Serfling1980}.

Let $\Delta:=\E D_{12}$ denote the population Gini mean difference. Since
$D_{12}^{2}=(X_1-X_2)^2$ and $X_1,X_2$ are i.i.d. with variance $\sigma^2$,
\[
\E D_{12}^{2}=\E(X_1-X_2)^{2}=2\sigma^{2}.
\]
Therefore
\[
\zeta_2=2\sigma^{2}-\Delta^{2},
\qquad
\zeta_1=\E\{\eta(X_1)^{2}\}-\Delta^{2}.
\]
In the notation of \citet[eqs.~(2.1)--(2.3) and~(4.3)]{ZengaPolisicchioGreselin2004},
their $\Delta$ is the same population Gini mean difference, their mean-deviation
function $D(x)=\E|x-X|$ is $\eta(x)$ here, and their expected squared
mean-deviation functional $\mathcal F=\E\{D(X)^2\}$ is
$\E\{\eta(X_1)^2\}$. Thus $\zeta_1=\mathcal F-\Delta^2$, and
\[
n\Var\{\GMD(X)\}\to 4\zeta_1
\qquad (n\to\infty).
\]
Equivalently, $4\zeta_1$ is the asymptotic variance of
$\sqrt n\{\GMD(X)-\Delta\}$
\citep[eqs.~(8)--(9) and the display following eq.~(9)]{GerstenbergerVogel2015}.
It is also the i.i.d. special case of the GMD long-run variance in
\citet[eq.~(3)]{GerstenbergerVogelWendler2020}, since their projection
$\varphi(x)=\E|x-Y|-g(F)$ is $\eta(x)-\Delta$ in the present notation.

\subsection*{Exact variances and adjacency fraction}
\begin{proposition}[Exact component variances]\label{prop:exact}
Let $X_1,\dots,X_n$ $(n\ge 2)$ be i.i.d. with finite variance, and let $\Pi$ be the independent uniformly random permutation of Section~\ref{sec:setup}. Then
\begin{align}
\Var\{\overline{\mathrm{MR}}(X,\Pi)\}
&=\frac{(n-1)\,\zeta_2+2(n-2)\,\zeta_1}{(n-1)^2},
\label{eq:varMR_exact}\\
\Var\{\GMD(X)\}
&=\frac{2\bigl\{\zeta_2+2(n-2)\,\zeta_1\bigr\}}{n(n-1)},
\label{eq:varGMD_exact}\\
\E\bigl[\Var\{\overline{\mathrm{MR}}(X,\Pi)\mid X\}\bigr]
&=\frac{n-2}{n(n-1)}\left\{\zeta_2-\frac{2(n-2)}{n-1}\,\zeta_1\right\}.
\label{eq:adjcomp_exact}
\end{align}
Dividing each expression by $d_2^{2}$ gives the corresponding components for $T=\overline{\mathrm{MR}}/d_2$.
\end{proposition}

Equation~\eqref{eq:varGMD_exact} is the exact finite-sample variance of a
degree-two $U$-statistic with kernel $h(x,y)=|x-y|$, written in the projection
components $(\zeta_1,\zeta_2)$ \citep{Hoeffding1948,Serfling1980}. For a coefficient check, put $A=\sum_{i<j}D_{ij}$. Then
\[
\Var(A)=\binom{n}{2}\zeta_2+6\binom{n}{3}\zeta_1,
\]
because there are $\binom{n}{2}$ variance terms and, within each triple of
observations, three unordered pairs of pairwise differences sharing one
observation; the factor 2 comes from the covariance terms in the variance
expansion. Dividing by $\binom{n}{2}^2$ gives \eqref{eq:varGMD_exact}. It was first
derived for the Gini mean difference by \citet{Nair1936} and rederived by
\citet{Lomnicki1952} (see \citealp{David1968} for the history). In the notation of
\citet[eqs.~(2.1)--(2.3), (3.2), and~(4.3)]{ZengaPolisicchioGreselin2004},
\[
\mathcal F=\zeta_1+\Delta^2,
\qquad
\sigma^2=\frac{\zeta_2+\Delta^2}{2}.
\]
Substituting these into their equation~(3.2) gives
\[
\frac{4}{n(n-1)}
\left\{
\frac{\zeta_2+\Delta^2}{2}
+(n-2)(\zeta_1+\Delta^2)
-\frac{2n-3}{2}\Delta^2
\right\}
=
\frac{2\{\zeta_2+2(n-2)\zeta_1\}}{n(n-1)},
\]
which is \eqref{eq:varGMD_exact}. In the notation of
\citet[eqs.~(8)--(9)]{GerstenbergerVogel2015}, write
$g_{\mathrm{GV}}=g(F)=\Delta$ and let $J$ denote their triple-integral
functional. Their equation~(8) becomes \eqref{eq:varGMD_exact} under
\[
\zeta_2=2\sigma^2-g_{\mathrm{GV}}^2,
\qquad
\zeta_1=\sigma^2+4J-g_{\mathrm{GV}}^2,
\]
because
\[
2\{\zeta_2+2(n-2)\zeta_1\}
=
4(n-1)\sigma^2+16(n-2)J-2(2n-3)g_{\mathrm{GV}}^2.
\]
Equation \eqref{eq:varMR_exact}, which is specific to successive absolute
differences, follows from the one-dependence calculation below, and
\eqref{eq:adjcomp_exact} is the difference between \eqref{eq:varMR_exact} and
\eqref{eq:varGMD_exact}. The
adjacency component vanishes at $n=2$, as it should, because both orderings of a
pair give the same moving range.

\begin{corollary}\label{cor:adjfrac}
Under the conditions of Proposition~\ref{prop:exact} with $\zeta_2>0$, the adjacency fraction \eqref{eq:adjfrac} (equivalently, its unscaled version for $\overline{\mathrm{MR}}$) depends on the sampling distribution only through $\rho_1$:
\begin{equation}
\mathrm{AdjFrac}(n)
=\frac{n-2}{n}\cdot\frac{(n-1)-2(n-2)\rho_1}{(n-1)+2(n-2)\rho_1},
\qquad
\lim_{n\to\infty}\mathrm{AdjFrac}(n)=\frac{1-2\rho_1}{1+2\rho_1}.
\label{eq:adjfrac_rho}
\end{equation}
\end{corollary}

Because \eqref{eq:adjcomp_exact} is an expected conditional variance, it is nonnegative for every $n$. For $n>2$ the prefactor in \eqref{eq:adjcomp_exact} is positive, so
\[
\zeta_2-\frac{2(n-2)}{n-1}\,\zeta_1\ge 0.
\]
Because this inequality holds for every $n>2$, letting $n\to\infty$ gives $\zeta_2\ge 2\zeta_1$. Since $\zeta_1\ge 0$ and $\zeta_2>0$ for a nondegenerate distribution, $\rho_1\in[0,1/2]$, and the limiting fraction in \eqref{eq:adjfrac_rho} lies in $[0,1]$.

Using $\zeta_2=2\sigma^2-\Delta^2$, $\zeta_1=\E\{\eta(X_1)^2\}-\Delta^2$, and standard Gini mean-difference moments, representative values are
\[
\begin{array}{c|cccc}
F & \mathrm{Uniform} & \mathrm{Normal} & t_5 & \mathrm{Exponential}\\
\hline
\rho_1 & 1/10 & 0.2239 & 0.3146 & 1/3\\
\mathrm{AdjFrac}(\infty) & 2/3 & 0.3813 & 0.2276 & 1/5
\end{array}
\]
The Uniform value follows from the Rectangular case of
\citet[Sec.~4.3]{ZengaPolisicchioGreselin2004}: after setting $b-a=1$,
$\Delta=1/3$, $\mathcal F=7/60$, and $\sigma^2=1/12$, so
$\zeta_1=1/180$, $\zeta_2=1/18$, and $\rho_1=1/10$. The Exponential value
follows from \citet[Sec.~4.2]{ZengaPolisicchioGreselin2004}: with rate parameter
$\theta$, $\Delta=1/\theta$, $\mathcal F=4/(3\theta^2)$, and
$\sigma^2=1/\theta^2$, so $\zeta_1=1/(3\theta^2)$,
$\zeta_2=1/\theta^2$, and $\rho_1=1/3$. The $t_5$ value uses
\citet[Theorem~2 and Tables~4--5]{GerstenbergerVogel2015}, with
$\zeta_1=\mathrm{ASV}(g_n)/4$ and $\zeta_2=2\sigma^2-g(F)^2$. The Normal value
follows from the constants in \eqref{eq:zetas_normal} below. Online Resource~1
confirms all four pairs by simulation. Online Resource~2 gives a complementary
nested-permutation check: for each reference distribution it repeatedly draws a
sample vector $X$, estimates
$\Var\{\overline{\mathrm{MR}}(X,\Pi)\mid X\}$ by randomizing permutations
conditional on that vector, and compares the resulting Monte Carlo averages with
\eqref{eq:varMR_exact}--\eqref{eq:adjfrac_rho}.

\subsection*{Proofs and Normal-reference check}
\paragraph{Proof of \eqref{eq:varMR_exact}.}
By exchangeability, $\Var\{\overline{\mathrm{MR}}(X,\Pi)\}$ equals the variance along any fixed ordering, so set $D_t:=D_{t,t+1}=|X_{t+1}-X_t|$ for $t=1,\dots,n-1$. The sequence $D_1,\dots,D_{n-1}$ is stationary and one-dependent, with $\Var(D_t)=\zeta_2$ and $\Cov(D_t,D_{t+1})=\zeta_1$: adjacent differences share the single observation $X_{t+1}$, while differences with disjoint index pairs are independent. The standard variance formula for the sum of a stationary one-dependent sequence \citep[e.g.,][]{BrockwellDavis1991} then gives
\[
\Var\Bigl(\sum_{t=1}^{n-1}D_t\Bigr)
=(n-1)\,\zeta_2+2(n-2)\,\zeta_1,
\]
and dividing by $(n-1)^2$ gives \eqref{eq:varMR_exact}. \hfill$\square$

\paragraph{Proof of \eqref{eq:adjcomp_exact} and Corollary~\ref{cor:adjfrac}.}
By Proposition~\ref{prop:gmd}, $\E\{\overline{\mathrm{MR}}(X,\Pi)\mid X\}=\GMD(X)$, so the law of total variance \eqref{eq:lotv} gives
\[
\E\bigl[\Var\{\overline{\mathrm{MR}}(X,\Pi)\mid X\}\bigr]
=\Var\{\overline{\mathrm{MR}}(X,\Pi)\}-\Var\{\GMD(X)\}.
\]
Over the common denominator $n(n-1)^2$, subtracting \eqref{eq:varGMD_exact} from \eqref{eq:varMR_exact} yields
\[
\E\bigl[\Var\{\overline{\mathrm{MR}}(X,\Pi)\mid X\}\bigr]
=\frac{(n-1)(n-2)\,\zeta_2-2(n-2)^2\,\zeta_1}{n(n-1)^{2}},
\]
which rearranges to \eqref{eq:adjcomp_exact}. Dividing this expression by \eqref{eq:varMR_exact} and simplifying yields \eqref{eq:adjfrac_rho}, proving Corollary~\ref{cor:adjfrac}; the limit follows by letting $n\to\infty$. \hfill$\square$

\paragraph{Normal-reference constants.}
Under $N(\mu,\sigma^{2})$ sampling, the bracket in \eqref{eq:VarTbar_normal_sigma} equals $a(n-2)+(\pi-2)$ with $a=\pi/3+2\sqrt{3}-4$, so
\[
\Var\{\GMD\}
=\frac{4}{\pi}\,\Var\{\bar T\}
=\frac{4\sigma^{2}}{\pi}\,\frac{a(n-2)+(\pi-2)}{n(n-1)}.
\]
Matching the constant and $(n-2)$ coefficients against \eqref{eq:varGMD_exact} gives
\begin{equation}
\zeta_2=\Bigl(2-\frac{4}{\pi}\Bigr)\sigma^2,
\qquad
\zeta_1=\Bigl(\frac{1}{3}+\frac{2\sqrt{3}-4}{\pi}\Bigr)\sigma^2.
\label{eq:zetas_normal}
\end{equation}
Substituting \eqref{eq:zetas_normal} into \eqref{eq:varMR_exact} and multiplying
by $1/d_2^{2}=\pi/4$ gives \eqref{eq:VarT_normal_sigma}; before this final
scaling, this is Hoel's unit-variance formula for $\Var(w)$
\citep[eq.~(10)]{Hoel1946}, with an additional factor $\sigma^2$ for general
scale. The value of $\zeta_1$ also agrees with Hoel's shared-difference moment
$\E\,D_{12}D_{23}=\tfrac13+2\sqrt{3}/\pi$ for $\sigma=1$
\citep[eq.~(8)]{Hoel1946}, since
$\zeta_1=\E\,D_{12}D_{23}-\Delta^{2}$ with $\Delta=2/\sqrt{\pi}$. Finally,
$\rho_1=\zeta_1/\zeta_2=a/(2\pi-4)\approx 0.2239$, and the limit in
\eqref{eq:adjfrac_rho} reproduces \eqref{eq:AdjFrac_limit_normal} exactly.

\bibliographystyle{plainnat}
\bibliography{refs_v14}

\end{document}